\documentclass[a4paper,10pt]{amsart}
\usepackage[margin=1.11in]{geometry}
\usepackage[utf8]{inputenc}
\usepackage{amsfonts,amssymb,amsmath,amstext,amscd,epsfig}
\usepackage{amsthm,cases,color,graphicx,graphics}
\usepackage{stackengine}
\usepackage{mathtools}
\usepackage[colorlinks=true,pdfborder={0 0 0}]{hyperref}
\hypersetup{urlcolor=blue,citecolor=teal,linkcolor=purple}
\numberwithin{equation}{section}
\newtheorem{lemma}{Lemma}[section]
\newtheorem{theorem}[lemma]{Theorem}

\newtheorem{proposition}[lemma]{Proposition}
\newtheorem{remark}{Remark}

\newcommand{\thishypname}{}
\newtheorem*{generichypothesis}{\thishypname}

\usepackage{tikz}
\def\asum{\mathop{\mathchoice
    {\stackon[-3.8ex]{\displaystyle\sum}
      {\smash{\rule{.4pt}{4ex}}}}
    {\stackon[-2.6ex]{\textstyle\sum}
      {\smash{\rule{.4pt}{2.9ex}}}}
    {\stackon[-1.9ex]{\scriptstyle\sum}
      {\smash{\rule{.4pt}{2.2ex}}}}
    {\stackon[-1.4ex]{\scriptscriptstyle\sum}
      {\smash{\rule{.4pt}{1.7ex}}}}
}}
\stackMath
\newcommand{\R}{\mathbb{R}}

\newcommand{\Z}{\mathbb{Z}}
\newcommand{\T}{\mathbb{T}}
\newcommand{\G}{\mathbb{G}}

\newcommand{\N}{\mathbb{N}}
\newcommand{\sX}{\mathsf{X}}

\newcommand{\dd}{\, {\rm d}}
\newcommand{\dt}{\frac{{\rm d}}{{\rm d}t}}

\newcommand{\cC}{\mathcal{C}}
\newcommand{\cR}{\mathcal{R}}
\newcommand{\cL}{\mathcal{L}}

\newcommand{\cM}{\mathcal{M}}
\newcommand{\cO}{\mathcal{O}}
\newcommand{\fB}{\mathfrak{B}}

\newcommand{\Lip}{\operatorname{Lip}}

\title[A consistency-stability approach to hydrodynamic limit\dots]{A
  consistence-stability approach to hydrodynamic limit of interacting
  particle systems on lattices}

\author{Angeliki Menegaki}
\address{Institut des hautes études Scientifiques, 35 Rte de
  Chartres, 91440, Bures-sur-Yvette, France}
   \email{Menegaki@ihes.fr}

\author{Clément Mouhot}
\address{DPMMS, University of Cambridge, Wilberforce Road,
  Cambridge CB3 0WA, UK}
   \email{C.Mouhot@dpmms.cam.ac.uk}

\begin{document}
\maketitle

\begin{abstract}
  This is a review based on the presentation done at the seminar
  Laurent Schwartz in December $2021$. It is announcing results in the
  forthcoming~\cite{MMM22}. This work presents a new simple
  quantitative method for proving the hydrodynamic limit of a class of
  interacting particle systems on lattices. We present here this
  method in a simplified setting, for the zero-range process and the
  Ginzburg-Landau process with Kawasaki dynamics, in the parabolic
  scaling and in dimension $1$. The rate of convergence is
  quantitative and uniform in time. The proof relies on a
  consistence-stability approach in Wasserstein distance, and it
  avoids the use of the ``block estimates''. 
\end{abstract}

\tableofcontents

\section{The general method}
\label{sec:intro}

We consider the hydrodynamic limit of interacting particle systems on
a lattice. The problem is to show that under an appropriate scaling of
time and space, the local particle densities of a stochastic lattice
gas converge to the solution of a macroscopic partial differential
equation.  We first present our method abstractly and then sketch
applications to two concrete models: the zero-range process (ZRP) and
the Ginzburg Landau process with Kawasaki dynamics (GLK). The
hydrodynamic limit is known at a qualitative level for all these
models under both hyperbolic and parabolic scalings for the ZRP and
under parabolic scaling for the GLK,
see~\cite{GPV88,Yau91,Reza,KL99}. However finding quantitative error
estimates had remained an important opened question, as well as
understanding the long-time behaviour of the hydrodynamic limit. First
results towards quantitative error, in the particular case of the
Ginzburg-Landau process with Kawasaki dynamics in dimension $1$, were
obtained in the two-parts work~\cite{DMOW18a,DMOW18b}, which builds
upon partial progresses in~\cite{GOVW09}.

\subsection{Set up and notation}

We denote by $\sX$ the state space at a given site (number of
particles, spin, etc.), which will in this paper be $\N$ (ZRP) or $\R$
(GLK). Consider the discrete torus $\T_N^d$ and the corresponding
phase space of particle configurations $\sX_N:=\sX^{\T^d
  _N}$. Variables in $\T^d_N$ are called \emph{microscopic} and
denoted by $x,y,z$, whereas variables in the limit continuous torus
$\T^d$ are called \emph{macroscopic} and denoted by $u$; finally
particle \emph{configurations} in $\sX_N$ are denoted by $\eta$. The
canonical embeddding $\T_N^d \to \T^d$, $x \mapsto \frac{x}{N}$ means
the macroscopic distance between sites of the lattice is
$\frac{1}{N}$. Given a particle configuration $\eta\in \sX_N$, we
define the \emph{empirical measure}
\begin{equation}
  \label{eq:emp_meas_dirac}
  \alpha_\eta^N := \asum_{x\in\T^d_N} \eta_x
  \delta_{\frac{x}{N}} \in \mathcal{M}_+(\T^d).
\end{equation}
where $\eta_x$ denotes the value of $\eta$ at $x \in \T^d_N$, and
$\mathcal{M}_+(\T^d)$ is the space of positive Radon measures on the
torus, and $\asum$ denotes the ``average sum'', here
$N^{-d} \sum_{x\in\T^d _N}$.

At the microscopic level, the interacting particle system evolves
through a stochastic process and the time-dependent probability
measure describing the law of $\eta$ is denoted by
$\mu_t^N \in P(\sX_N)$.  We consider a linear operator
$\cL_N : C_b(\sX_N) \rightarrow C_b(\sX_N)$ generating uniquely a
Feller semigroup $e^{t \cL_N}$ on $P(\sX_N)$ (see~\cite[Chapter
1]{Liggett1985}) so that given
$\mu_0^N \in P(\sX_N)$ the solution
$\mu_t^N = e^{t \cL_N} \mu^N_0 \in P(\sX_N)$ satisfies
\begin{equation}
  \label{eq:particle_evol}
  \forall \, \Phi \in C_b(\sX_N), \quad
  \dt\langle \Phi, \mu^N_t \rangle = \langle \cL_N \Phi, \mu^N_t
  \rangle,
\end{equation}
where $C_b(\sX_N)$ denotes continuous bounded real-valued functions
and $\langle \cdot, \cdot \rangle$ denotes the duality bracket between
$C_b(\sX_N)$ and $P(\sX_N)$.

At the macroscopic level, we consider a map 
$\cL_\infty : \cM_+(\G_\infty) \rightarrow \cM_+(\G_\infty)$ (in
general unbounded and nonlinear) and
the evolution problem
\begin{equation}
  \label{eq:limitPDE1} 
  \partial_t f_t = \cL_\infty f_t, \quad f_{t =0}=f_0.
\end{equation}

A measure $\mu^N \in P(\sX_N)$ is called \emph{invariant}
for~\eqref{eq:particle_evol} if
\begin{equation*}
  \forall \, \Phi \in C_b(\sX_N), \quad
  \big\langle \mu^N, \mathcal{L}_N \Phi \big\rangle = 0.
\end{equation*}

We also denote $\Lip(\sX_N)$ the Lipschitz functions
$\Phi : \sX_N \to \R$ with respect to the (normalised) $\ell^1$ norm:
for every $\eta, \zeta \in \sX_N$,
$|\Phi(\eta)- \Phi(\zeta)| \leq C_\Phi \asum_{x\in \T_N^d}
|\eta_x-\zeta_x|$, and we denote the smallest such constant $C_\Phi$
by $[\Phi]_{\Lip(\sX_N)} \in \R_+$.

\subsection{Abstract assumptions}

We make the following assumptions
on~\eqref{eq:particle_evol}-\eqref{eq:limitPDE1}: \smallskip

\noindent
\textbf{(H0) Local equilibrium structure.} There are
$n_\lambda: \text{Conv}(\sX) \to \R_+$ depending on $\lambda \in \R$
($\text{Conv}$ denotes the convex hull) and
$\sigma : \text{Conv}(\sX) \to \R$ so that: (i)
$n_\lambda ^{\otimes \T^d_N}$ is invariant on $\sX_N$ for each
$\lambda$, and (ii) for any $\rho \in \text{Conv}(\sX)$,
$\mathbb E_{n_{\sigma(\rho)}}[\eta_x] = \rho$. We then define, given
a macroscopic profile $f$ on $\T^d$, the \emph{local Gibbs measure}
\begin{equation*}
  \vartheta^N_f(\eta) := \nu_{\sigma \left( f \left(
        \frac{\cdot}{N} \right) \right)} ^N(\eta) \quad
  \text{where} \quad \nu_{F}^N(\eta) := \prod_{x \in \G_N}
  n_{F(x)}(\eta(x)).
\end{equation*}

\begin{figure}[!ht]
  \includegraphics[scale=0.8]{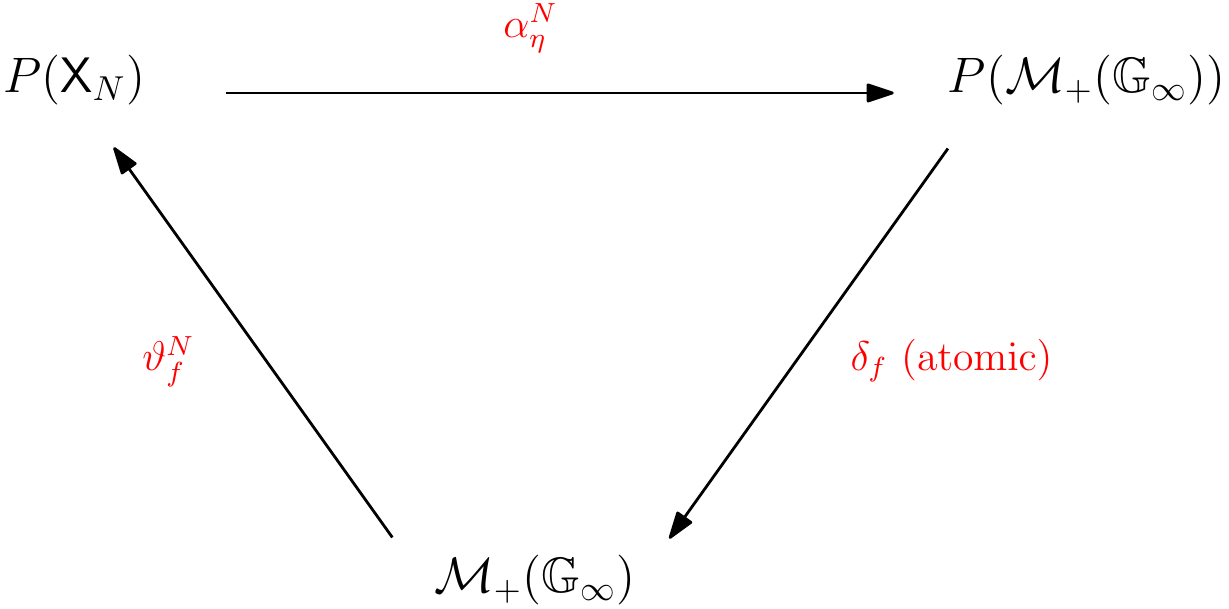}
  \caption{The functional setting.}
  \label{fig:setting}
\end{figure}

The two maps $\eta \mapsto \alpha_\eta ^N$ and
$f \mapsto \vartheta^N_f$ allow comparisons between the
microscopic and macroscopic scales, as summarized in
Figure~\ref{fig:setting}.
\smallskip

\noindent
\textbf{(H1) Microscopic stability.}  The semigroup
$e^{t\mathcal{L}_N}$ satisfies
\begin{align}
  \forall \, \Phi \in \Lip(\sX_N), \quad \left[ e^{t\mathcal{L}_N}
  \Phi \right]_{\Lip(\sX_N) } \leq  \left[  \Phi
  \right]_{\Lip(\sX_N)}.
\end{align}  

\noindent
\textbf{(H2) Macroscopic stability.} There is a Banach space
$\fB \subset \cM_+(\G_\infty)$ so that~\eqref{eq:limitPDE1} is locally
well-posed in $\fB$; given the maximal time of existence
$T_m \in (0,+\infty]$ we denote for $t \in [0,T_m)$,
$R(t) := \left\| f_t - f_{\infty} \right\|_\fB$
when~\eqref{eq:limitPDE1} has a unique stationary solution
$f_\infty \in \fB$ with mass $\int_{\T^d} f_\infty = \int_{\T^d} f_0$,
otherwise we denote $R(t) := \left\| f_t \right\|_\fB$.  \smallskip

\noindent
\textbf{(H3) Consistency.}  There is a \emph{consistency error}
$\epsilon(N) \to 0$ as $N \to \infty$ so that for $T \in [0,T_m)$
\begin{equation*}
  \frac1T \int_0 ^T
  \int_0 ^t \left\langle \left( e^{(t-s)\cL_N} \Phi \right), \left[
      \cL_N^* \left( \frac{{\rm d} \vartheta_{f_s} ^N}{{\rm d}
          \nu^N_\infty}\right) - \frac{{\rm d}}{{\rm d}s} \left(
        \frac{{\rm d} \vartheta_{f_s} ^N}{{\rm d} \nu^N_\infty}
      \right) \right] {\rm d} \nu^N_\infty \right\rangle \dd s \dd t 
  \le \epsilon(N) [\Phi]_{\Lip(\sX_N)} \int_0 ^t R(s) \dd s
\end{equation*}
for any $\Phi \in \Lip(\sX_N)$, where $\nu^N_\infty$ is an equilibrium
measure. 

\subsection{The abstract strategy}

\begin{theorem}
  \label{thm:main}
  Consider~\eqref{eq:particle_evol}-\eqref{eq:limitPDE1} with the
  assumptions {\bf (H0)--(H1)--(H2)--(H3)}. Let
  $\phi \in C^\infty (\mathbb{T}^d)$, $\mu_0^N \in P_1(X_N)$ for all
  $N\geq 1$, $f_0 \in \mathcal{B}$. Then
  \begin{equation}
    \label{eq:cvg}
    \forall \, T \in [0,T_m), \quad
    \frac1T \int_0 ^T \left\| \mu_t ^N - \vartheta_{f_t}^N
    \right\|_{\Lip^*} \dd t \lesssim  \epsilon(N) \int_0^T R(s)
    \dd s + \left\| \mu_0 ^N- \vartheta_{f_0} ^N \right\|_{\Lip^*}.
  \end{equation}
\end{theorem}

\begin{remark}
  Note that $\| \mu_t^N - \vartheta_{f_t}^N \|_{\Lip^*} \to 0$ as
  $N \to \infty$ implies that the empirical
  measure~\eqref{eq:emp_meas_dirac} sampled from the law $\mu^N_t$
  satisfies
  \begin{equation}
    \label{eq:hydro_profile}
    \forall \, \phi \in C_b(\G), \ \forall \, \epsilon >0, \
    \forall \, t \ge 0,\quad
    \lim_{N\rightarrow\infty} \mu^N_t\left( \left\{ |\langle
        \alpha^N_\eta, \varphi \rangle - \langle f_t, \varphi \rangle
        | > \epsilon \right\} \right) = 0
  \end{equation}
  with a rate of convergence (thus recovering quantitatively results
  from~\cite{GPV88}):
  \begin{align*}
    & \mu^N_t\left( \left\{ |\langle
      \alpha^N_\eta, \varphi \rangle - \langle f_t, \varphi \rangle
      | > \epsilon \right\} \right) \\ 
    & \le \mu^N_t\left( \left\{ \langle
      \alpha^N_\eta, \varphi \rangle \ge \langle f_t, \varphi \rangle
      + \epsilon \right\} \right) + \mu^N_t\left( \left\{ \langle
      \alpha^N_\eta, \varphi \rangle \le \langle f_t, \varphi \rangle
      - \epsilon \right\} \right) \\
    & \le \int_{\sX_N} \left[ F_\epsilon ^+ \left( \left\langle \phi,
      \alpha_\eta ^N \right\rangle \right)  -
      F_\epsilon ^+  \left( \left\langle \phi, f_t \right\rangle
      \right) \right] \dd \mu_t^N + \int_{\sX_N} \left[ F_\epsilon ^-
      \left( \left\langle \phi, \alpha_\eta ^N
      \right\rangle \right)  - F^-_\epsilon  \left( \left\langle \phi,
      f_t \right\rangle \right) \right] \dd \mu_t^N 
  \end{align*}
  where $F_\epsilon ^\pm$ are mollified version of the
  characteristic functions of respectively
  $\{ z \ge \left\langle \phi, f_t \right\rangle + \epsilon \}$
  and
  $\{ z \le \left\langle \phi, f_t \right\rangle - \epsilon
  \}$, which yields
  \begin{equation*}
    \sup_{t \in [0,T]} \mu^N_t\left( \left\{ |\langle
        \alpha^N_\eta, \varphi \rangle - \langle f_t, \varphi \rangle
        | > \epsilon \right\} \right) \lesssim \epsilon^{-1}
    \| \mu_t^N - \vartheta_{f_t}^N \|_{\Lip^*} + \epsilon^{-2} N^{-d}.
  \end{equation*}
\end{remark}

\section{Concrete applications}

We apply the abstract result to two archetypical models, the
zero-range process (ZRP), and the Ginzburg-Landau process with
Kawasaki dynamics (GLK).

\subsection{The ZRP} In this case, the state space at each site
is $\sX=\N$. Given the choice of a \emph{transition function}
$p \in P(\T^d_N)$ with $p(0)=0$ and a \emph{jump rate function}
$g : \N \rightarrow \R_+$, the base generator $\hat \cL_N$ writes
\begin{equation}
  \label{eq:zrp-gen}
  \forall \, \Phi \in C_b(\sX_N), \ \forall \, \eta \in \sX_N, \quad
  \hat \cL_N \Phi(\eta) := \sum_{x,y \in \T^d_N} p(y-x) g(\eta_x)
  \left[\Phi(\eta^{xy}) - \Phi(\eta) \right]
\end{equation}
where $\eta^{xy}$ is defined as before. The local equilibrium
structure of {\bf (H0)} is given by
\begin{align}
  \label{eq:zrp-gibbs}
  & n_\lambda(k):= \frac{\lambda^k}{g(k)!
    Z(\lambda)} \quad \text{ with } \quad
    Z(\lambda) := \sum_{k=0} ^{+\infty} \frac{\lambda^k}{g(k)!}\\
  & \sigma \ \text{ is defined implicitely by } \ \sigma(\rho)
    \frac{Z'(\sigma(\rho))}{Z(\sigma(\rho))} \equiv \rho
\end{align}
denoting $g(k)!:=g(k) g(k-1) \cdots g(1)$. The pair $(g,\sigma)$ thus
constructed satisfies
$\mathbb E_{n_{\sigma(\alpha)}}[g]=\sigma(\alpha)$. When
$f \equiv \rho \in [0,+\infty)$ is constant, the local Gibbs measure
$\vartheta_{\rho}^N=\nu_{\sigma(\rho)}^N$ is invariant with average
number of particles $\rho$. The \emph{mean transition rate} is defined
by $\gamma := \sum_{x \in \Z^d} x p(x) \in \R^d$. When
$\gamma \not =0$, the first non-zero asymptotic dynamics as
$N \to \infty$ is given by the hyperbolic scaling
$\cL_N := N \hat \cL_N$, and the corresponding expected limit equation
is $\partial_t f = \gamma \cdot \nabla [\sigma(f)]$. When $\gamma =0$,
the first non-zero asymptotic dynamics as $N \to \infty$ is the given
by the parabolic scaling $\cL_N := N^2 \hat \cL_N$, and the
corresponding limit equation is formally
\begin{equation}
  \label{eq:zrp-diff}
  \partial_t f = \Delta_a [\sigma(f)] \quad \text{ with } \quad
  \Delta_a := \sum_{i,j=1} ^d a_{ij} \partial^2_{ij}
  \quad \text{and} \quad
  a_{ij} := \sum_{x \in \Z^d} p(x) x_i x_j.
\end{equation}

We make the following assumptions on the jump rate function
$g:\N \rightarrow [0,\infty)$.
\smallskip

\noindent {\bf (HZRP)} The jump rate $g$ satisfies
$g(0) = 0$, $g(n)>0$ for all $n>0$, is non-decreasing, uniformly
Lipschitz $\sup_{n \ge 0} |g(n+1)-g(n)|<+\infty$, and there are
$n_0 > 0$ and $\beta > 0$ such that $g(n')-g(n) \ge \beta$ for
any $n'\ge n+n_0$.
\medskip

The main result on the ZRP is:
\begin{theorem}[Hydrodynamic limit for the ZRP]
  \label{theo:ZRP}
  Consider $\hat \cL_N$ defined in~\eqref{eq:zrp-gen} with $g$
  satisfying {\bf (HZRP)}. Let $d=1$, $f_0 \in C^3(\T)$ with
  $f_0 \ge \delta>0$, and $\mu_0^N \in P_1(\sX_N)$ for all $N \ge
  1$. Assume $\gamma=0$, define $\mu^N_t = e^{t N^2 \hat \cL_N}$ and
  $f_t \in C([0,T),C^3(\T^d))$ solution to~\eqref{eq:zrp-diff}, then
  the following convergence holds (with quantitative constants)
  \begin{equation}
    \label{eq:dec-zrp}
    \sup_{T \ge 0} \frac1T \int_0 ^T \left\| \mu_t ^N - \vartheta_{f_t}^N
    \right\|_{\Lip^*} \dd t \lesssim N^{-\frac18} + \left\| \mu_0 ^N-
      \vartheta_{f_0} ^N \right\|_{\Lip^*}.
  \end{equation}
\end{theorem}

\subsection{The GLK} In this case, the state space at each site
is $\sX=\R$. Given the choice of a \emph{single-site potential}
$V \in C^2(\R)$, the base generator $\hat \cL_N$ writes
\begin{equation}
  \label{eq:glk-gen}
  \hat \cL_N \Phi(\eta) :=
  \frac{1}{2}\sum_{x \sim y \in \T^d_N} \left( \frac{\partial}{\partial
      \eta_x}  - \frac{\partial}{\partial \eta_y} \right)^2
  -\frac{1}{2} \sum_{x \sim y \in \T^d_N} \left[ V'(\eta_x) - V'(\eta_y)
  \right] \left( \frac{\partial}{\partial \eta_x} -
    \frac{\partial}{\partial \eta_y} \right)
\end{equation}
where $x \sim y$ denotes neighbouring sites. The local
equilibrium structure is given by
\begin{align*}
  & n_\lambda(r):= \frac{e^{\lambda r-V(r)}}{Z(\lambda)}
  \quad \text{ with } \quad Z(\lambda) := \int_\R e^{\lambda r -
    V(r)} \dd r \\
  & \sigma \ \text{ is defined implicitely by } \
    \frac{Z'(\sigma(\rho))}{Z(\sigma(\rho))} \equiv \rho.
\end{align*}
When $f \equiv \rho \in \R$ is constant, the local Gibbs measure
$\vartheta_{\rho}^N=\nu_{\sigma(\rho)}^N$ is invariant with average
spin $\rho$. The hyperbolic scaling formally leads to zero and the
parabolic scaling $\cL_N := N^2 \hat \cL_N$ formally leads to
\begin{equation}
  \label{eq:glk-diff}
  \partial_t f = 2 \Delta[\sigma(f)].
\end{equation}

We assume that the single-site potential satisfies
\smallskip

\noindent {\bf (HGLK)} The potential $V$ is $C^2$ and decomposes as
$V(u)=V_0(u)+V_1(u)$ with $V_0''(u) \ge \kappa$ for all $u \in \R$ for
some $\kappa>0$ and $\| V_1\|_{W^{1,\infty}(\R)} \lesssim 1$.
\medskip

This assumption is similar with those in~\cite{GOVW09, DMOW18a,
  Fat13}. One can take for example a double-well potential,
provided it is uniformly convex at infinity.

\begin{theorem}[Hydrodynamic limit for the GLK]
  \label{theo:GLK}
  Consider $\cL_N$ defined in~\eqref{eq:glk-gen} with $V$ satisfying
  {\bf (HGLK)}. Let $d=1$, $f_0 \in C^3(\T^d)$ and
  $\mu_0^N \in P_1(\sX_N)$ for all $N \ge 1$. Define
  $\mu^N_t = e^{t N^2 \hat \cL_N}$ and
  $f_t \in C([0,+\infty),C^3(\T^d))$ the global solution
  to~\eqref{eq:glk-diff}, then the following convergence holds (with
  quantitative constants) 
  \begin{equation}
    \label{eq:dec-glk}
    \sup_{T \ge 0} \frac1T \int_0 ^T \left\| \mu_t ^N - \vartheta_{f_t}^N
    \right\|_{\Lip^*} \dd t \lesssim N^{-\frac18} + \left\| \mu_0 ^N-
      \vartheta_{f_0} ^N \right\|_{\Lip^*}.
  \end{equation}
 \end{theorem}

\section{The abstract strategy}
\label{sec:abstract}

In this section we sketch the proof of Theorem~\ref{thm:main}. Let
$f_t$ be a solution to~\eqref{eq:limitPDE1}. Given $0<\ell<N$, we
denote by $\eta^{\ell}$ for the local $\ell$-average
$\eta^{\ell}_x := \asum_{|y-x|\leq \ell} \eta_y$.

Denote by $F^N_t := \frac{{\rm d} \mu^N_t}{{\rm d}\nu^N_\infty}$ and
$G^N_t := \frac{{\rm d} \vartheta^N_{f_t}}{{\rm d}\nu^N_\infty}$ the
densities with respect to $\nu^N_\infty$, and write
\begin{align*}
  \frac{ \dd}{ \dd t} \Big( F^N_t - G_t ^N \Big) =
  \mathcal{L}_N^* \Big( F^N_t - G_t ^N \Big) + \left(
  \mathcal{L}_N ^* G_t ^N - \partial_t G_t ^N \right)
\end{align*}
so that Duhamel's formula yields
\begin{align*}
  F^N_t - G_t ^N = e^{t \mathcal{L}_N ^*} \left( F^N_0 - G_0 ^N
  \right) + \int_0^t e^{(t-s)\mathcal{L}_N ^*} \left( \mathcal{L}_N ^*
  G_s ^N - \partial_s G_s ^N \right) \dd s.
\end{align*}
Take $\Phi \in \Lip(\sX_N)$ with
$ \| \Phi \|_{\operatorname{Lip}(\sX_N)} \leq 1$ and integrate
the above equation to get
\begin{align*}
  & \int_{\sX_N} \Phi  \Big( F^N_t - G_t ^N \Big) \dd \nu^N_\infty \\
  & = \underbrace{\int_{\sX_N} \left( e^{t \mathcal{L}_N} \Phi \right)
    \Big( F^N_0 - G_0 ^N \Big) \dd \nu^N_\infty}_{I_1(t)} +
    \underbrace{\int_{\sX_N} \int_0^t \left( e^{(t-s) \mathcal{L}_N}
    \Phi  \right) \left( \mathcal{L}_N G_s ^N - \partial_s G_s  ^N
    \right) \dd \nu^N_\infty \dd s}_{I_2(t)}.
\end{align*}
\textbf{(H1)} implies
$I_1(t) \lesssim \| \mu^N_0 - \vartheta_{f_0} ^N \|_{\Lip^*}$ and
\textbf{(H3)} implies
$\frac1T \int_0 ^T I_2(t) \dd t \le \epsilon(N) \int_0 ^T R(s) \dd s$,
which implies the conclusion of Theorem~\ref{thm:main}.

\section{Proof for the ZRP}
\label{sec:sep}

In this section we prove Theorem~\ref{theo:ZRP}) (hydrodynamical limit
for the ZRP). Note for this model $\cL_N=\cL_N^*$ is symmetric with
respect to equilibrium measures.

Given $f_t \in C^3(\T^d)$ with $f> \delta$, $\delta>0$, and
$\rho := \int_{\T^d} f$, the density of the local Gibbs measure
relatively to the invariant measure with mass $\rho$ is:
\begin{equation}
  \label{eq:zrp-psi}
  G_t^N (\eta) := \frac{{\rm d} \vartheta_f^N(\eta)}{{\rm d}
    \vartheta_{\rho}^N(\eta)} = \prod_{x \in \T^d_N}
  \left(\frac{\sigma \left(f_t \left( \frac{x}{N} \right)
      \right)}{\sigma(\rho)} \right)^{\eta(x)}
  \left(\frac{ Z(\sigma \left(f_t \left( \frac{x}{N} \right)
      \right))}{Z(\sigma(\rho))} \right)^{-1}.
\end{equation}
where the function $\sigma(r)$ is defined by
$\langle n_{\sigma(r)}, \eta(x) \rangle = r$ and the
\emph{partition function}
$Z:[0,\lambda^\ast) \rightarrow \mathbb{R}$ is defined
in~\eqref{eq:zrp-gibbs}, with $\lambda^* \in [0,+\infty]$
denoting the radius of convergence of the series.

It is proved in~\cite[Chapter~2, Section~3]{KL99} that assumption~{\bf
  (HZRP)} on $g$ implies that
$\sigma = R^{-1} : [0,\infty) \rightarrow [0,\infty)$ is well-defined
and strictly increasing, with
\begin{equation*}
  R(\lambda) = \lambda \partial_\lambda \log(Z(\lambda)) =
  \frac{1}{Z(\lambda)} \sum_{n\geq0} \frac{n\lambda^n}{g(n)!}.
\end{equation*}
Then the building block $n_\rho$ of the Gibbs measure satisfies
$\langle n_{\sigma(\rho)}, g(\eta(x)) \rangle = \sigma(\rho)$.
Moreover {\bf (HZRP)} implies that the function $\sigma$ is $C^\infty$
with uniform bound on all derivatives on $\R_+$, with Lipschitz
constant less than $g^\ast$, and with
$\inf_{\lambda >0} \lambda^{-1} \sigma(\lambda) >0$ (in particular
$\sigma'(0)>0$), and that the invariant measure has exponential moment
bounds,
see~\cite[Corollary~3.6]{KL99}.

\subsection{Microscopic Stability -- {\bf (H1)}}
\label{subsec:microscopic stability_ZRP}

We use the coupling of~\cite{Liggett1985,Reza}. Let
\begin{equation}
  \label{def:sep-coupling-op}
  \begin{split}
    \widetilde{\mathcal{L}}_N\Psi(\eta, \zeta) :=& \sum_{x,y \in
      \T^d_N} p(y-x) \Big( g(\eta_x) \wedge g(\zeta_x) \Big) \Big[
    \Psi(\eta^{xy}, \zeta^{xy}) -
    \Psi(\eta,\zeta) \Big] \\
    & + \sum_{x,y \in \T^d_N} p(y-x) \Big(g(\eta_x) - g(\eta_x) \wedge
    g(\zeta_x) \Big)
    \Big[ \Psi(\eta^{xy}, \zeta) - \Psi(\eta, \zeta) \Big] \\
    & + \sum_{x,y \in \T^d_N} p(y-x) \Big( g(\zeta_x) - g(\eta_x)
    \wedge g(\zeta_x) \Big) \Big[ \Psi(\eta,\zeta^{xy}) - \Psi(\eta,
    \zeta) \Big].
  \end{split}
\end{equation}
for a two-variable test function $\Psi(\eta,\zeta)$. Then
$\widetilde{\mathcal{L}}_N\Phi(\eta) = \hat \cL_N \Phi(\eta)$ and
$\widetilde{\mathcal{L}}_N\Phi(\zeta) = \hat \cL_N \Phi(\zeta)$, and
{\bf (H1)} follows from the fact that $e^{t\widetilde{\mathcal{L}}_N}$
preserves sign and the inequality
\begin{equation*}
  \widetilde{\mathcal{L}}_N \left( \sum_{z \in \T^d_N} |\eta_z
    - \zeta_z| \right) \le 0.
\end{equation*}
To prove the latter inequality, we compute
\begin{align*} 
  \widetilde{\mathcal{L}}_N \left( \sum_{z \in \T^d_N} |\eta_z
  - \zeta_z| \right)
  & = \sum_{x,y \in \T^d_N} p(y-x) \Big( g(\eta_x) - g(\eta_x) \wedge
    g(\zeta_x) \Big) \\
  & \hspace{2cm} \times \Big[ |\eta^{xy}_x- \zeta_x| +
    |\eta^{xy}_y- \zeta_y| - |\eta_x- \zeta_x| -
    |\eta_y-\zeta_y| \Big] \\
  & \quad + \sum_{x,y \in \T^d_N} p(y-x) \Big( g(\zeta_x) - g(\eta_x)
    \wedge g(\zeta_x) \Big) \\
  & \hspace{2cm} \times \Big[ |\eta_x- \zeta^{xy}_x| +
    |\eta_y - \zeta^{xy}_y| - |\eta_x - \zeta_x| -
    |\eta_y -\zeta_y| \Big].
\end{align*}
When $g(\eta_x) - g(\eta_x) \wedge g(\zeta_x) >0$ necessarily
$\eta_x- \zeta_y \ge 1$ and
\begin{equation*}
  \Big[ |\eta^{xy}_x- \zeta_x| + |\eta^{xy}_y- \zeta_y| -
  |\eta_x- \zeta_x| - |\eta_y-\zeta_y| \Big] \le 0.
\end{equation*}
When $g(\zeta_x) - g(\eta_x) \wedge g(\zeta_x) >0$ necessarily
$\zeta_x - \eta_x \ge 1$ and
\begin{equation*}
  \Big[ |\eta_x- \zeta^{xy}_x| + |\eta_y- \zeta^{xy}_y| -
  |\eta_x- \zeta_x| - |\eta_y-\zeta_y| \Big] \le 0.
\end{equation*}

\subsection{Macroscopic stability -- {\bf (H2)}}

In the parabolic scaling the limit PDE is the nonlinear diffusion
equation~\eqref{eq:zrp-diff}. We take $\mathcal{B} = C^3$ with its
standard infinity Banach norm. The proof that this norm remains
uniformly bounded in time is classical in dimension $d=1$ (using the
bounds on $\sigma$), and $f_t \in [\delta, 1-\delta]$ for all times by
maximal principle. Moreover $f_t \to \rho$ exponentially fast as $t
\to \infty$ in $\mathcal{B}$.

\subsection{Consistency estimate -- {\bf (H3)}} 

Let $\gamma = 0$ and the dimension $d=1$.
\begin{proposition}
  \label{prop:zrp-consist-tspt}
  Given the solution $f_t \in C^3(\T^d)$ to~\eqref{eq:zrp-diff} with
  $f \ge \delta$, $\delta>0$, and $\rho := \int_{\T^d} f$, and
  $G_t ^N$ defined in~\eqref{eq:zrp-psi}, we have for every
  $\Phi \in \operatorname{Lip}(\sX_N)$ with
  $[\Phi]_{\operatorname{Lip}(\sX_N)} \le 1$
  \begin{align*}
    \frac1T \int_0 ^T I_t ^N \dd t := \frac1T \int_0 ^T \int_0 ^t
    \left\langle \left( e^{(t-s)\cL_N} \Phi \right), \left[ \cL_N
    G^N_s  - \frac{{\rm d}}{{\rm d}s} G^N_s \right] {\rm d}
    \nu^N_\infty \right\rangle \dd s \dd t =  \cO \left(
    N^{-\frac18}\right)  
  \end{align*}
  where the constant depends on the estimates in {\bf (H2)}.
\end{proposition}
 
\begin{proof}
  We start by computing
  \begin{equation*}
    \cL_N G^N_s  - \frac{{\rm d}}{{\rm d}s} G^N_s = \sum_{x \in
      \T^d_N} A_x ^N G^N_s
  \end{equation*}
  with (note that $f_t \to \rho$ exponentially fast)
  \begin{align*}
    A_x ^N
    & := N^2 \sum_{y \in \T^d_N} p(y-x) g(\eta_x) \left(
      \frac{\sigma\left( f_t\left(\frac{y}{N}\right)
      \right)}{\sigma\left( f_t\left(\frac{x}{N}\right)
      \right)} -1 \right) - \eta_x \frac{\sigma'\left(
      f_t\left(\frac{x}{N}\right) \right)}{\sigma\left(
      f_t\left(\frac{x}{N}\right) \right)}
      \Delta_a[\sigma(f)]\left(\frac{x}{N}\right) \\
    & = \frac{g(\eta_x)}{\sigma\left( f_t\left(\frac{x}{N}\right)
      \right)} \Delta_a[\sigma(f)]\left(\frac{x}{N}\right) - \eta_x
      \frac{\sigma'\left( f_t\left(\frac{x}{N}\right)
      \right)}{\sigma\left( f_t\left(\frac{x}{N}\right) \right)}
      \Delta_a[\sigma(f)]\left(\frac{x}{N}\right) +
      \cO\left(\frac{e^{-Cs}}{N}\right)
  \end{align*}
  for some $C>0$. Since (conservation of mass)
  \begin{equation*}
    \int_{\sX_N} \left( \sum_{x \in \T^d_N} A_x ^N G^N_s \right) \dd
    \nu_\infty ^N = \int_{\sX_N} \left( \sum_{x \in \T^d_N} A_x ^N \right) \dd
    \vartheta_{f_s} ^N =0,
  \end{equation*}
  we can replace $\Phi_{t-s} := e^{(t-s)\cL_N} \Phi$ by
  \begin{equation*}
    \tilde \Phi_{t,s} := e^{(t-s)\cL_N} \Phi - {\bf
      E}_{\vartheta_{f_s}^N}[e^{(t-s)\cL_N} \Phi]
  \end{equation*}
  and use the Lipschitz bound on $e^{(t-s)\cL_N} \Phi$ (microscopic
  stability) to get
   \begin{equation*}
      I_t ^N = \int_0 ^t \int_{\sX_N} \tilde \Phi_{t,s}(\eta) \left(
        \sum_{x \in \T^d_N} \tilde A^N_x \right) \dd \vartheta_{f_s}
      ^N + \cO\left(\frac{1}{N}\right)
   \end{equation*}
   with $\tilde A^N_x$ defined by (note that it has zero average
   against ${\rm d} \vartheta_{f_s} ^N$)
   \begin{equation*}
     \tilde A^N_x :=  \left\{ g(\eta_x) - \sigma\left(
         f_t\left(\frac{x}{N}\right) \right) - \sigma'\left(
         f_t\left(\frac{x}{N}\right) \right) \left[ \eta_x -
         f_t\left(\frac{x}{N}\right) \right] \right\}
     \frac{\Delta_a[\sigma(f)]\left(\frac{x}{N}\right)}{\sigma\left(
         f_t\left(\frac{x}{N}\right) \right)}.
   \end{equation*}
   We then form sub-sum over non-overlapping cubes of size
   $\ell \in \{1,\dots,N\}$ (this intermediate scale factor $\ell$
   will be chosen later in terms of $N$). Let $\cR^d_N \subset \T^d_N$
   be a net of centers of non-overlapping cubes of the form
   $\cC_x := \{ y \in \T^d_N \ : \ \| x - y \|_\infty \le \ell
   \}$. Then
   \begin{align*}
     I_t ^N
     & = \sum_{x \in \cR_N^d}  \int_0 ^t \int_{\sX_N} \tilde
       \Phi_{t,s}(\eta) \left( \sum_{y \in \cC_x} \tilde A^N_y \right)
       \dd \vartheta_{f_s} ^N + \cO\left(\frac{1}{N}\right) \\
     & = (2\ell+1)^d \sum_{x \in \cR_N^d}  \int_0 ^t \int_{\sX_N}
       \tilde \Phi_{t,s}(\eta) \hat A^N_x \dd \vartheta_{f_s} ^N +
       \cO\left(\frac{1}{N}\right)
   \end{align*}
   with the $\hat A^N_x$ defined by 
   \begin{equation*}
     \hat A^N_x :=  \left\{ \langle g(\eta) \rangle_{\cC_x} -
       \sigma\left( f_t\left(\frac{x}{N}\right) \right) -
       \sigma'\left( f_t\left(\frac{x}{N}\right) \right) \left[
         \langle \eta \rangle_{\cC_x} - f_t\left(\frac{x}{N}\right)
       \right] \right\}
     \frac{\Delta_a[\sigma(f)]\left(\frac{x}{N}\right)}{\sigma\left(
         f_t\left(\frac{x}{N}\right) \right)}
   \end{equation*}
   where $\langle F(\eta) \rangle_{\cC_x}$, for $F=F(\eta_x)$, denotes
   taking the average over the cube $\cC_x$. Note that the average of
   $\hat A^N _x$ against ${\rm d} \vartheta_{f_s} ^N$ is
   $\cO(e^{-Cs}\ell/N)$. Then
   \begin{align*}
     & \sum_{x \in \cR_N^d}  \int_0 ^t \int_{\sX_N} \tilde \Phi_{t,s}
       \hat A^N_x \dd \vartheta_{f_s} ^N \\
     & = \sum_{x \in \cR_N^d}  \int_0 ^t \int_{\sX_N} \left( \tilde
       \Phi_{t,s} - \Pi_x ^N \tilde \Phi_{t,s} \right) \hat A^N_x \dd
       \vartheta_{f_s} ^N + \sum_{x \in \cR_N^d}  \int_0 ^t
       \int_{\sX_N} \Pi_x ^N \tilde \Phi_{t,s} \hat A^N_x \dd
       \vartheta_{f_s} ^N \\
     & = \sum_{x \in \cR_N^d}  \int_0 ^t \int_{\sX_N} \left(
       \Phi_{t-s} - \Pi_x ^N \Phi_{t-s} \right) \hat A^N_x \dd
       \vartheta_{f_s} ^N \\
     & \hspace{2.5cm} + \sum_{x \in \cR_N^d}  \int_0 ^t \int_{\sX_N}
       \left( \Pi_x ^N \Phi_{t-s} - {\bf E}_{\vartheta_{f_s} ^N
       }[\Pi_x ^N \Phi_{t-s}] \right) \hat A^N_x \dd \vartheta_{f_s}
       ^N =: J^N_t + \tilde J^N_t
   \end{align*}
   where $\Pi_x ^N$ projects on the average over configurations
   $\Omega_m := \left\{ \tilde \eta \, : \, \langle \tilde \eta
     \rangle_{\cC_x} = m \right\}$ with the same mass in the cube
   $\cC_x$ (and does not touch the other sites):
   \begin{equation}
     \label{eq:pi-average}
       \Pi_x ^N \varphi(\eta) = [\Pi_x ^N \varphi](\langle \eta
       \rangle_{\cC_x} ) = \int_{\Omega_{\langle \eta
           \rangle_{\cC_x}}} \varphi(\tilde \eta) \dd
       \nu^{\ell,\langle \eta
         \rangle_{\cC_x}}(\tilde \eta)
   \end{equation}
   for a function $\varphi$ on $\sX^{\cC_x}$. To estimate the first
   term $J^N_t$ we first approximate the measure $\vartheta_{f_s} ^N$
   on $\cC_x$ by the equilibrium measure with local mass $f_t(x/N)$,
   and denote it by $\bar \vartheta_{f_s}$ (note that the
   approximation is made differently for each cube and depends on $x$,
   even if it is written explicitly). This produces an error
   $\cO(\ell^{d+1}/N)$ (using the Lipschitz regularity of $\Phi_{t-s}$
   and the exponential convergence $f_t \to \rho$ to get uniform in
   time bounds). We then apply the Poincaré
   inequality~\cite[Theorem~1.1]{MR1415232} in the cube $\cC_x$ (whose
   constant is independent of the number of particles and proportional
   to the size of the cube) and the law of large number
   $\| \hat A^N_x \|_{L^2(\bar \vartheta^N_{f_s})} = \cO(e^{-Cs}
   \ell^{-d/2})$ (using uniform bounds on the second moment of
   $\bar \vartheta_{f_s}$):
   \begin{align*}
     J^N_t
     & \le \sum_{x \in \cR_N^d}  \int_0 ^t \| \Phi_{t-s} - \Pi_x ^N
       \Phi_{t-s} \|_{L^2(\bar \vartheta_{f_s} ^N)} \| \hat A^N_x
       \|_{L^2(\bar \vartheta_{f_s} ^N)} \dd s + \cO\left(
       \frac{\ell^{d+1}}{N} \right) \\
     & \lesssim \ell^{1-\frac{d}2} \sum_{x \in \cR_N^d}  \int_0 ^t
       \sqrt{\bar D^\ell_x\left(\Phi_{t-s}\right)} e^{-Cs} \dd s + \cO\left(
       \frac{\ell^{d+1}}{N} \right) \\
     & \lesssim \ell^{1-\frac{d}2} N^{\frac{d}{2}} \int_0 ^t \left(
       \sum_{x \in \cR_N^d} \bar D^\ell_x \left(\Phi_{t-s}\right)
       \right)^{\frac12} e^{-Cs} \dd s + \cO\left( \frac{\ell^{d+1}}{N}
       \right)
   \end{align*}
   where $\bar D^\ell_x(\Phi)$ is the Dirichlet form on the cube
   $\cC_x$ with respect to the measure $\bar \vartheta_{f_s} ^N$:
   \begin{equation*}
     \bar D^\ell_x(\Phi) := \sum_{y,z \in \cC_x} \int_{\sX_N} p(z-y)
     g(\eta_y) \left[ \Phi(\eta^{yz}) - \Phi(\eta) \right]^2 \dd \bar
     \vartheta_{f_s} ^N.
   \end{equation*}
   Then we change back the measure $\bar \vartheta_{f_s} ^N$ in each
   box, which produces (using the Lipschitz regularity of $\Phi$) an
   error $\ell^{3/2} N^{-1/2}$), and we compute
   \begin{equation*}
     \frac{1}{2N^2} \dt \int_{\sX_N} \Phi_{t-s}(\eta)^2 \dd 
     \vartheta_{f_s} ^N \le - \sum_{x \in \cR_N^d} D^\ell_x
     \left(\Phi_{t-s}\right) + \cO\left( \frac{1}{N^2} \right)
   \end{equation*}
   (with $D^\ell$ denoting the Dirichlet form for
   $\vartheta_{f_s} ^N$), where the last error accounts for the small
   default of self-adjointness. We deduce (in dimension $d=1$) that
   \begin{equation*}
     \int_0 ^T J_t ^N \dd t \lesssim T^{\frac12} \left( \frac{\ell}{N}
     \right)^{1-\frac{d}2}  + \cO\left( \frac{T \ell^{d+1}}{N} \right).
   \end{equation*}

   To control the second term $\tilde J_t ^N$, we first use the
   \emph{equivalence of ensemble} in~\cite[Appendix~II,
   Corollary~1.7]{KL99} on the local equilibrium measure
   $\bar \vartheta_{f_s} ^N$ (using uniform exponential moment
   bounds):
   \begin{equation}
     \label{eq:equiv-ens}
     \langle g(\eta) \rangle_{\cC_x} = \sigma \left( \langle \eta
       \rangle_{\cC_x} \right) + \cO\left( \frac{1}{\ell^d} \right).
   \end{equation}
   Second we remark that the Lipschitz regularity of $\Phi_{t-s}$
   implies that
   $\Pi_x ^N \Phi_{t-s} - {\bf E}_{\vartheta_{f_s} ^N }[\Pi_x ^N
   \Phi_{t-s}] = \cO(\ell^d N^{-d})$, and since the average of
   $\hat A_x ^N$ with respect to $\vartheta_{f_s} ^N$ is
   $\cO(\ell/N)$, we can write
   \begin{equation*}
     \tilde J_t ^N = \sum_{x \in \cR_N^d}  \int_0 ^t \int_{\sX_N} \left(
       \Pi_x ^N \Phi_{t-s}[\langle \eta \rangle_{\cC_x}] - \Pi_x ^N
       \Phi_{t-s}\left[f_s\left( \frac{x}{N} \right)\right] \right)
     \hat A^N_x \dd \vartheta_{f_s} ^N + \cO\left( \frac{\ell}{N}
     \right).
   \end{equation*}
   Third, we remark that the Lipschitz regularity of $\Phi_{t-s}$
   (with constant $N^{-d}$) implies a Lipschitz regularity of its
   averaged projection $\Pi_x ^N \Phi_{t-s}$ with constant
   $\ell^d N^{-d}$, with respect to the local mass. Indeed, given
   $0=m \le m' <+\infty$, pick any pair of configuration
   $(\eta_0,\zeta_0)$ with $\langle \eta_0 \rangle_{\cC_x} =m$,
   $\langle \zeta_0 \rangle_{\cC_x} =m'$ and $\eta_0 \le \zeta_0$
   (such configuration trivially exists since $m \le m'$). Then we
   consider the initial coupling $\delta_{(\eta_0,\zeta_0)}$ on
   $\Omega_m \times \Omega_{m'}$ which has $\ell^1$ cost $m'-m$. Then
   we evolve it along the flow of the coupling operator
   $e^{t \tilde \cL_N}$. The marginals respectively converge to
   $\nu^{\ell,m}$ and $\nu^{\ell,m'}$ (convergence to equilibrium of
   the oiriginal evolution). Since the evolution by the coupling
   operator does not increase the Wasserstein distance, we deduce
   $W_1(\nu^{\ell,m},\nu^{\ell,m'})\le m'-m$. An optimal coupling
   $\Pi$ associated to this distance thus satisfies
   \begin{equation*}
     m'-m \le \int_{\Omega_m \times \Omega_{m'}} \left( \asum_{x \in \T^d_N}
       |\eta_x - \zeta_x| \right) \Pi(\eta,\zeta) \le m'-m
   \end{equation*}
   where the first inequality follows from Jensen's inequality. Thus
   the Jensen's inequality is saturated which implies that the cost
   does not change sign on the support of $\Pi$, i.e. $\eta \le \zeta$
   in the support. We then compute
   \begin{align*}
     \Pi_x ^N \Phi_{t-s}(m') - \Pi_x ^N \Phi_{t-s}(m) 
     & = \int_{\Omega_{m'}} \Phi_{t-s}(\zeta) \dd \nu^{\ell,m'}(\zeta)
       - \int_{\Omega_m} \Phi_{t-s}(\eta) \dd \nu^{\ell,m}(\eta) \\
     & = \int_{\Omega_m \times \Omega_{m'}} \left[ \Phi_{t-s}(\zeta) -
       \Phi(\eta) \right]  \dd \Pi(\eta,\zeta) 
   \end{align*}
   and since $\eta \le \zeta$ on the support of $\Pi$, $\|
   \zeta-\eta\|_{\ell^1(\cC_x)} = (m'-m) \ell^d$ and 
   \begin{equation*}
     \left| \Pi_x ^N \Phi_{t-s}(m') - \Pi_x ^N \Phi_{t-s}(m) \right|
     \le \frac{\ell^d}{N^d} |m'-m|.
   \end{equation*}
   We deduce (using~\eqref{eq:equiv-ens})
   \begin{multline*}
     \tilde J_t ^N \lesssim \frac{\ell^d}{N^d} \sum_{x \in \cR_N^d}  \int_0
     ^t \int_{\sX_N} \left| \langle \eta \rangle_{\cC_x} - f_s\left(
         \frac{x}{N} \right) \right| \times \\ \left| \sigma(\langle
       \eta \rangle_{\cC_x}) - \sigma\left( f_s\left( \frac{x}{N}
         \right) \right) - \sigma'\left( f_s\left( \frac{x}{N} \right)
       \right) \left[ \langle \eta \rangle_{\cC_x} - f_s\left(
           \frac{x}{N} \right) \right] \right| \dd \vartheta_{f_s} ^N
     e^{-Cs} \dd s
     \\ + \frac{1}{N^d} \sum_{x \in \cR_N^d}  \int_0 ^t \int_{\sX_N}
     \left| \langle \eta \rangle_{\cC_x} - f_s\left( \frac{x}{N}
       \right) \right| \dd \vartheta^N_{f_s} e^{-Cs} \dd s + \cO\left(
       \frac{\ell}{N} \right)
   \end{multline*}
   which yields by Taylor formula, the approximation of
   $\vartheta^N_{f_s}$ by $\bar \vartheta^N_{f_s}$, and the law of
   large numbers
   \begin{align*}
     \tilde J_t ^N
     & \lesssim \frac{\ell^d}{N^d} \sum_{x \in \cR_N^d}  \int_0 ^t
       \int_{\sX_N} \left| \langle \eta \rangle_{\cC_x} - f_s\left(
       \frac{x}{N} \right) \right|^3  \dd \vartheta_{f_s} ^N e^{-Cs}
       \dd s \\
     & \hspace{2cm} + \frac{1}{N^d} \sum_{x \in \cR_N^d}  \int_0 ^t
       \int_{\sX_N} \left| \langle \eta \rangle_{\cC_x} - f_s\left(
       \frac{x}{N} \right) \right| \dd \vartheta^N_{f_s} e^{-Cs} \dd s
       + \cO\left( \frac{\ell}{N} \right) e^{-Cs} \dd s \\
     & \lesssim \frac{\ell^d}{N^d} \sum_{x \in \cR_N^d}  \int_0 ^t
       \int_{\sX_N} \left| \langle \eta \rangle_{\cC_x} - f_s\left(
       \frac{x}{N} \right) \right|^3  \dd \bar \vartheta_{f_s} ^N
       e^{-Cs} \dd s\\
     & \hspace{2cm} + \frac{1}{N^d} \sum_{x \in \cR_N^d}  \int_0 ^t
       \int_{\sX_N} \left| \langle \eta \rangle_{\cC_x} - f_s\left(
       \frac{x}{N} \right) \right| \dd \bar \vartheta^N_{f_s} e^{-Cs}
       \dd s + \cO\left( \frac{\ell}{N} \right) \\
     & \lesssim \cO\left( \ell^{-\frac{3d}{2}} \right) +
       \cO\left( \frac{\ell}{N} \right).
   \end{align*}
   Combining all estimates we get (optimizing $\ell:=N^{1/4}$)
   \begin{equation}
     \label{eq:conclusion-proof}
     \frac1T \int_0 ^T I_t ^N \dd t \lesssim \left( \frac{1}{N} +
       \frac{\ell^{1+\frac{d}{2}}}{N^{1-\frac{d}{2}}} +
       \frac{\ell^{1+2d}}{N} + \frac{1}{\ell^{\frac{d}{2}}} +
       \frac{\ell}{N} \right) \lesssim \frac{1}{N^{\frac{1}{8}}}.
   \end{equation}
\end{proof}

\section{Proof for the GLK}
\label{sec:glk}

In this section we prove Theorem~\ref{theo:GLK} (hydrodynamic limit
for the GLK). Note again that for this model $\cL_N=\cL_N^*$ is
symmetric with respect to equilibrium measures. Given
$f_t \in C^3(\T^d)$ and $\rho := \int_{\T^d} f \in \R$, the density of
the local Gibbs measure relatively to the invariant measure with mass
$\rho$ is:
\begin{equation}
  \label{eq:glk-psi}
  G_t^N (\eta) := \frac{{\rm d} \vartheta_f^N(\eta)}{{\rm d}
    \vartheta_{\rho}^N(\eta)} = \prod_{x \in \T^d_N}
  e^{\left[\sigma \left(f_t \left( \frac{x}{N} \right) \right) -
        \sigma(\rho)\right] \eta_x}
    \frac{Z(\sigma(\rho))}{Z\left(\sigma\left(f_t \left( \frac{x}{N}
          \right)\right)\right)}.
\end{equation}
where the function $\sigma(r)$ is defined by
$\langle n_{\sigma(r)}, \eta_x \rangle = r$ and the \emph{partition
  function} $Z(\lambda) = \int_\R e^{\lambda r - V(r)} \dd r$ is
defined on $\R$. The uniform convexity of $V$ at infinity easily
implies bounds on some exponential moments of the invariant measure
and {\bf (HGLK)} implies that there exists $C>0$ so that
$ 0 < \frac{1}{C} \leq \sigma' \leq C < \infty $ (see~\cite[Lemma
41]{GOVW09} and~\cite[Lemma~5.1]{DMOW18a}).

\subsection{Microscopic stability -- {\bf (H1)}}

We define a ``coupling generator''
$\widetilde{\mathcal{L}}_N:C_b(\sX_N^2) \to C_b(\sX_N^2)$ by
\begin{equation}
  \label{eq:coupled_gen_GL}
  \begin{split}  
    \widetilde{\mathcal{L}}_N \Psi(\eta, \zeta) := 
    \sum_{x \sim y} \bigg( & \left[ \left(
        \frac{\partial}{\partial \eta_x} -
        \frac{\partial}{\partial \eta_y} \right)^* \left(
        \frac{\partial}{\partial \eta_x} -
        \frac{\partial}{\partial \eta_y} \right) \otimes 1
    \right] \Psi(\eta, \zeta) \\
    & + \left[ 1 \otimes \left( \frac{\partial}{\partial
          \zeta_x} - \frac{\partial}{\partial \zeta_y}
      \right)^* \left( \frac{\partial}{\partial \zeta_x} -
        \frac{\partial}{\partial \zeta_y} \right) \right]
    \Psi(\eta, \zeta) \\
    & + K \left( \frac{\partial}{\partial \eta_x} -
      \frac{\partial}{\partial \eta_y} \right) \otimes \left(
      \frac{\partial}{\partial \zeta_x} -
      \frac{\partial}{\partial \zeta_y} \right) \Psi(\eta,\zeta)
    \bigg)
  \end{split}
\end{equation} 
where $K>0$ is a constant to be chosen later and the adjoint is
taken in $L^2({\rm d}\vartheta^N_\rho)$ so 
\begin{align*}
  \hat \cL_N
  & = \sum_{x \sim y} \left( \frac{\partial}{\partial \eta_x}
    - \frac{\partial}{\partial \eta_y} \right)^2 - \left(
    V'(\eta_x) - V'(\eta_y) \right) \left(
    \frac{\partial}{\partial \eta_x} - \frac{\partial}{\partial
      \eta_y} \right) \\
  & = \sum_{x \sim y} \left( \frac{\partial}{\partial \eta_x} -
    \frac{\partial}{\partial \eta_y} \right)^*\left(
    \frac{\partial}{\partial \eta_x} - \frac{\partial}{\partial
      \eta_y} \right).
\end{align*}

Then for any $p \in (1,2]$ there is $K=K(p)>0$ (depending on $p$) so
that
\begin{align*}
  \widetilde{\mathcal{L}}_N\left( \sum_{x \in \T^d_N} |\eta_x- \zeta_x|^p
  \right)
  & = 2p(p-1)(2+4d) \sum_{x \in \T^d_N} |\eta_x - \zeta_x|^{p-2} \\ 
  & \qquad - 2 (p-1) \sum_{x \sim y} \left[ V_0'(\eta_x) -
    V_0'(\zeta_x) \right] (\eta_x-\zeta_x) |\eta_x - \zeta_x|^{p-1} \\
  & \qquad - 2 (p-1) \sum_{x \sim y} \left[ V_1'(\eta_x) -
    V_1'(\zeta_x) \right] (\eta_x-\zeta_x) |\eta_x - \zeta_x|^{p-1} \\
  & \qquad + K p(p-1) (2+4d) \sum_{x \in \T^d_N} |\eta_x -
    \zeta_x|^{p-2} \leq 0
\end{align*}
by using the assumptions on the potential: $V_0$ uniformly strictly
convex and $V_1 \in W^{1,\infty}$. This implies the weak contraction
of the evolution in $W_p$ ($p$-Wasserstein distance) for any
$p \in (1,2]$, and thus by limit in $W_1$. By duality this implies
that the evolution is weakly contractive for the dual Lipschitz norm.

\subsection{Macroscopic stability - {\bf (H2)}} 

The limit equation is similar to that of the ZRP and {\bf (H2)} is
proved in the same way.

\subsection{Consistency estimate - {\bf (H3)}} 
Let the dimension $d=1$.

\begin{proposition}
  \label{prop:glk-consist-tspt}
  Given the solution $f_t \in C^3(\T^d)$ to~\eqref{eq:zrp-diff}, and
  $\rho := \int_{\T^d} f$, and $G_t ^N$ defined in~\eqref{eq:glk-psi},
  we have for every $\Phi \in \operatorname{Lip}(\sX_N)$ with
  $[\Phi]_{\operatorname{Lip}(\sX_N)} \le 1$
  \begin{align*}
    \frac1T \int_0 ^T I_t ^N \dd t := \frac1T \int_0 ^T \int_0 ^t
    \left\langle \left( e^{(t-s)\cL_N} \Phi \right), \left[ \cL_N
    G^N_s  - \frac{{\rm d}}{{\rm d}s} G^N_s \right] {\rm d}
    \nu^N_\infty \right\rangle \dd s \dd t =  \cO \left(
    N^{-\frac18}\right)  
  \end{align*}
  where the constant depends on the estimates in {\bf (H2)}.
\end{proposition}

\begin{proof}
  The proof follows the same structure as for the ZRP. We start by
  computing
  \begin{equation*}
    \cL_N G^N_s  - \frac{{\rm d}}{{\rm d}s} G^N_s = \sum_{x \in
      \T^d_N} A_x ^N G^N_s
  \end{equation*}
  with (note again that $f_t \to \rho$ exponentially fast)
  \begin{align*}
    A_x ^N
    & =  \frac{N^2}{2} \sum_{y \sim x}  \Bigg[
      2 \sigma\left(f_s\left(\frac{x}{N}\right)\right) \left(  
      \sigma \left(f_s\left(\frac{x}{N}\right)\right) -
      \sigma \left(f_s\left(\frac{y}{N}\right)\right)  \right) \\
    & \hspace{3cm}
      - 2 V'(\eta_x)   \sigma \left(f_s\left(\frac{x}{N}\right)\right) -
      \sigma \left(f_s\left(\frac{y}{N}\right)\right)
      \Bigg]  \\
    & \hspace{5cm} - \sum_x \left( \eta_x -
      f_s\left(\frac{x}{N}\right) \right) \sigma'\left(
      f_s\left(\frac{x}{N}\right)
      \right) \Delta[\sigma(f)]\left(\frac{x}{N}\right)
    \\
    & = \Delta[\sigma(f)]\left(\frac{x}{N}\right)
      \Bigg[  V'(\eta_x) - 
      \sigma\left(f_s\left(\frac{x}{N}\right)\right)  -
      \sigma'\left(f_s\left(\frac{x}{N}\right)\right)
      \left( \eta_x- f_s\left(\frac{x}{N}\right)\right) \Bigg]
      + \cO\left(\frac{e^{-Cs}}{N}\right)
  \end{align*}
  for some $C>0$. By conservation of mass we replace again
  $\Phi_{t-s} := e^{(t-s)\cL_N} \Phi$ by
  \begin{equation*}
    \tilde \Phi_{t,s} := e^{(t-s)\cL_N} \Phi - {\bf
      E}_{\vartheta_{f_s}^N}[e^{(t-s)\cL_N} \Phi]
  \end{equation*}
  and use the Lipschitz bound {\bf (H1)} on $e^{(t-s)\cL_N} \Phi$ to
  get
   \begin{equation*}
      I_t ^N = \int_0 ^t \int_{\sX_N} \tilde \Phi_{t,s}(\eta) \left(
        \sum_{x \in \T^d_N} \tilde A^N_x \right) \dd \vartheta_{f_s}
      ^N + \cO\left(\frac{1}{N}\right)
   \end{equation*}
   with $\tilde A^N_x$ defined by (note that it has zero average
   against ${\rm d} \vartheta_{f_s} ^N$)
   \begin{equation*}
     \tilde A^N_x :=  \Delta[\sigma(f)]\left(\frac{x}{N}\right)
      \left[ V'(\eta_x) - \sigma\left(f\left(\frac{x}{N}\right) \right)
      - \sigma'\left(f_s\left(\frac{x}{N}\right)\right) \left( \eta_x
        - f\left(\frac{x}{N}\right)\right) \right].
   \end{equation*}
   We again form sub-sum over non-overlapping cubes of size
   $\ell \in \{1,\dots,N\}$, with $\cR^d_N \subset \T^d_N$ a net of
   centers of cubes
   $\cC_x := \{ y \in \T^d_N \ : \ \| x - y \|_\infty \le \ell
   \}$. Then
   \begin{align*}
     I_t ^N
     & = \sum_{x \in \cR_N^d}  \int_0 ^t \int_{\sX_N} \tilde
       \Phi_{t,s}(\eta) \left( \sum_{y \in \cC_x} \tilde A^N_y \right)
       \dd \vartheta_{f_s} ^N + \cO\left(\frac{1}{N}\right) \\
     & = (2\ell+1)^d \sum_{x \in \cR_N^d}  \int_0 ^t \int_{\sX_N}
       \tilde \Phi_{t,s}(\eta) \hat A^N_x \dd \vartheta_{f_s} ^N +
       \cO\left(\frac{1}{N}\right)
   \end{align*}
   with the $\hat A^N_x$ defined by (and
   $\langle F(\eta) \rangle_{\cC_x}$ again denotes the average over
   the cube $\cC_x$)
   \begin{equation*}
     \hat A^N_x :=  \Delta[\sigma(f)]\left(\frac{x}{N}\right)
      \left[ \langle V'(\eta) \rangle_{\cC_x} -
        \sigma\left(f\left(\frac{x}{N}\right) \right)
      - \sigma'\left(f_s\left(\frac{x}{N}\right)\right) \left( \langle
        \eta \rangle_{\cC_x} - f\left(\frac{x}{N}\right)\right) \right].
   \end{equation*}
   (Note again that the average of $\hat A^N _x$ against
   ${\rm d} \vartheta_{f_s} ^N$ is $\cO(e^{-Cs}\ell/N)$.) Then
   \begin{align*}
     & \sum_{x \in \cR_N^d}  \int_0 ^t \int_{\sX_N} \tilde \Phi_{t,s}
       \hat A^N_x \dd \vartheta_{f_s} ^N \\
     & \qquad = \sum_{x \in \cR_N^d}  \int_0 ^t \int_{\sX_N} \left( \tilde
       \Phi_{t,s} - \Pi_x ^N \tilde \Phi_{t,s} \right) \hat A^N_x \dd
       \vartheta_{f_s} ^N + \sum_{x \in \cR_N^d}  \int_0 ^t
       \int_{\sX_N} \Pi_x ^N \tilde \Phi_{t,s} \hat A^N_x \dd
       \vartheta_{f_s} ^N \\
     & \qquad = \sum_{x \in \cR_N^d}  \int_0 ^t \int_{\sX_N} \left(
       \Phi_{t-s} - \Pi_x ^N \Phi_{t-s} \right) \hat A^N_x \dd
       \vartheta_{f_s} ^N \\
     & \hspace{3.5cm} + \sum_{x \in \cR_N^d}  \int_0 ^t \int_{\sX_N}
       \left( \Pi_x ^N \Phi_{t-s} - {\bf E}_{\vartheta_{f_s} ^N
       }[\Pi_x ^N \Phi_{t-s}] \right) \hat A^N_x \dd \vartheta_{f_s}
       ^N =: J^N_t + \tilde J^N_t
   \end{align*}
   where $\Pi_x ^N$ again averages over $\Omega_m$ (and does not touch
   the other site) as in~\eqref{eq:pi-average}.

   To estimate the first term $J^N_t$ we again approximate the measure
   $\vartheta_{f_s} ^N$ on $\cC_x$ by the equilibrium measure with
   local mass $f_t(x/N)$, and denote it by $\bar \vartheta_{f_s}$
   (note that the approximation is made differently for each cube and
   depends on $x$, even if it is written explicitly). This produces an
   error $\cO(\ell^{d+1}/N)$ (using the Lipschitz regularity of
   $\Phi_{t-s}$ and the exponential convergence $f_t \to \rho$ to get
   uniform in time bounds). We then apply the Poincaré
   inequality~\cite[Theorem~2]{MR1233852} in the cube $\cC_x$ (whose
   constant is independent of the number of particles and proportional
   to the size of the cube) and the law of large number
   $\| \hat A^N_x \|_{L^2(\bar \vartheta^N_{f_s})} = \cO(e^{-Cs}
   \ell^{-d/2})$ (using uniform bounds on the second moment of $\bar
   \vartheta^N_{f_s}$)
   \begin{align*}
     J^N_t
     & \le \sum_{x \in \cR_N^d}  \int_0 ^t \| \Phi_{t-s} - \Pi_x ^N
       \Phi_{t-s} \|_{L^2(\bar \vartheta_{f_s} ^N)} \| \hat A^N_x
       \|_{L^2(\bar \vartheta_{f_s} ^N)} \dd s + \cO\left(
       \frac{\ell^{d+1}}{N} \right) \\
     & \lesssim \ell^{1-\frac{d}2} \sum_{x \in \cR_N^d}  \int_0 ^t
       \sqrt{\bar D^\ell_x\left(\Phi_{t-s}\right)} e^{-Cs} \dd s + \cO\left(
       \frac{\ell^{d+1}}{N} \right) \\
     & \lesssim \ell^{1-\frac{d}2} N^{\frac{d}{2}} \int_0 ^t \left(
       \sum_{x \in \cR_N^d} \bar D^\ell_x \left(\Phi_{t-s}\right)
       \right)^{\frac12} e^{-Cs} \dd s + \cO\left( \frac{\ell^{d+1}}{N}
       \right)
   \end{align*}
   where $\bar D^\ell_x(\Phi)$ is the Dirichlet form on the cube
   $\cC_x$ with respect to the measure $\bar \vartheta_{f_s} ^N$:
   \begin{equation*}
     \bar D^\ell_x(\Phi) := \sum_{y \sim z \in \cC_x} \int_{\sX_N} 
     \left[ \partial_{\eta_x} \Phi(\eta) - \partial_{\eta_y}
       \Phi(\eta) \right]^2 \dd \bar \vartheta_{f_s} ^N.
   \end{equation*}
   Then we use the entropy production
   \begin{equation*}
     \frac{1}{2N^2} \dt \int_{\sX_N} \Phi_{t-s}(\eta)^2 \dd 
     \vartheta_{f_s} ^N \le - \sum_{x \in \cR_N^d} D^\ell_x
     \left(\Phi_{t-s}\right) + \cO\left( \frac{1}{N^2} \right)
   \end{equation*}
   as before to deduce that
   \begin{equation*}
     \int_0 ^T J_t ^N \dd t \lesssim T^{\frac12} \left( \frac{\ell}{N}
     \right)^{1-\frac{d}2}  + \cO\left( \frac{T \ell^{d+1}}{N} \right).
   \end{equation*}

   To control the second term $\tilde J_t ^N$, we first use the
   \emph{equivalence of ensemble} in~\cite[Corollary~5.3]{MR1931585}
   on the local equilibrium measure $\bar \vartheta_{f_s} ^N$ (using
   bounds on some exponential moments):
   \begin{equation}
     \label{eq:equiv-ens}
     \langle V'(\eta) \rangle_{\cC_x} = \sigma \left( \langle \eta
       \rangle_{\cC_x} \right) + \cO\left( \frac{1}{\ell^d} \right).
   \end{equation}
   Second we remark that the Lipschitz regularity of $\Phi_{t-s}$
   implies that
   $\Pi_x ^N \Phi_{t-s} - {\bf E}_{\vartheta_{f_s} ^N }[\Pi_x ^N
   \Phi_{t-s}] = \cO(\ell^d N^{-d})$, and since the average of
   $\hat A_x ^N$ with respect to $\vartheta_{f_s} ^N$ is
   $\cO(\ell/N)$, we can write
   \begin{equation*}
     \tilde J_t ^N = \sum_{x \in \cR_N^d}  \int_0 ^t \int_{\sX_N} \left(
       \Pi_x ^N \Phi_{t-s}[\langle \eta \rangle_{\cC_x}] - \Pi_x ^N
       \Phi_{t-s}\left[f_s\left( \frac{x}{N} \right)\right] \right)
     \hat A^N_x \dd \vartheta_{f_s} ^N + \cO\left( \frac{\ell}{N}
     \right).
   \end{equation*}

   Third, we prove again that the Lipschitz regularity of $\Phi_{t-s}$
   (with constant $N^{-d}$) implies a Lipschitz regularity of its
   averaged projection $\Pi_x ^N \Phi_{t-s}$ with constant
   $\ell^d N^{-d}$, with respect to the local mass. Indeed, given
   $0=m < m' <+\infty$, pick any pair of configuration
   $(\eta_0,\zeta_0)$ with $\langle \eta_0 \rangle_{\cC_x} =m$,
   $\langle \zeta_0 \rangle_{\cC_x} =m'$ and $\eta_0 < \zeta_0$ (such
   configuration trivially exists since $m < m'$). Then consider the
   coupling on $\Omega_m \times \Omega_{m'}$ given by a product of
   smooth probability distributions localised around respectively
   $\delta_{\eta_0}$ and $\delta_{\zeta_0}$, so that the support of
   this coupling only contains strictly ordered $\eta < \zeta$. Then
   we evolve it along the flow of the coupling operator
   $e^{t \tilde \cL_N}$. The marginals respectively converge to
   $\nu^{\ell,m}$ and $\nu^{\ell,m'}$ (convergence to equilibrium of
   the oiriginal evolution). Arguing as for the ZRP, we deduce that
   $W_1(\nu^{\ell,m},\nu^{\ell,m'}) = m'-m$, and a corresponding
   optimal coupling $\Pi$ associated to this distance is so that the
   cost does not change sign on its support, i.e.  $\eta \le \zeta$ in
   the support. We deduce as for the ZRP that $\Pi_x ^N \Phi_{t-s}$ is
   $\ell^d N^{-d}$-Lipschitz.

   We finally deduce from~\eqref{eq:equiv-ens}, the Taylor formula,
   the approximation of $\vartheta^N_{f_s}$ by
   $\bar \vartheta^N_{f_s}$, and the law of large numbers, the same
   estimate on $\tilde J^N_t$ as for the ZRP, and finally the same
   estimate~\eqref{eq:conclusion-proof} on $I^N_t$, which concludes
   the proof.
\end{proof}

\bibliographystyle{alpha}
\bibliography{bibliography}

\bigskip

\end{document}